\documentclass[final,dvipdfm, 11pt]{amsart}
\usepackage{ amsmath, amsthm, amsfonts, graphicx}
\usepackage[dvipsnames,usenames]{color}
\usepackage{verbatim}
\usepackage{lineno}

\theoremstyle{plain}
\newtheorem{theorem}{Theorem}[section]
\newtheorem*{theorem*}{Theorem}

\newtheorem{pro}[theorem]{Proposition}
\newtheorem{Def}[theorem]{Definition}
\newtheorem{lem}[theorem]{Lemma}

\theoremstyle{definition}
\newtheorem*{Def*}{Definition}

\theoremstyle{plain}

\newtheorem{Rem}[theorem]{Remark}

\numberwithin{equation}{section}

\newcommand{\bpo}{\begin{pro}}
\newcommand{\epo}{\end{pro}}
\newcommand{\be}{\begin{equation}}
\newcommand{\ene}{\end{equation}}
\newcommand{\br}{\begin{Rem}}
\newcommand{\er}{\end{Rem}}
\newcommand{\bl}{\begin{lem}}
\newcommand{\el}{\end{lem}}
\newcommand{\bd}{\begin{Def}}
\newcommand{\ed}{\end{Def}}
\newcommand{\ben}{\begin{enumerate}}
\newcommand{\een}{\end{enumerate}}
\newcommand{\bp}{\begin{proof}}
\newcommand{\ep}{\end{proof}}
\newcommand{\beq}{\begin{equation*}}
\newcommand{\eeq}{\end{equation*}}
\newcommand{\bear}{\begin{eqnarray*}}
\newcommand{\eear}{\end{eqnarray*}}
\newcommand{\bt}{\begin{theorem}}
\newcommand{\et}{\end{theorem}}
\newcommand{\bst}{\begin{split}}
\newcommand{\est}{\end{split}}

\newcommand{\bal}{\begin{aligned}}
\newcommand{\eal}{\end{aligned}}

\renewcommand{\P}{\partial}
\newcommand{\F}[2]{\frac{#1}{#2}}
\newcommand{\la}{\langle}
\newcommand{\ra}{\rangle}

\newcommand{\R}{\mathbb{R}}

\newcommand{\bnb}{\bar{\nabla}}
\newcommand{\nb}{\nabla}

\newcommand{\RM}{Riemannian manifold}

\newcommand{\Lp}{Laplacian}

\newcommand{\Cd}{covariant derivative}

\newcommand{\Wp}{warped product manifold}

\newcommand{\ga}{\gamma}
\newcommand{\Ta}{\Theta}

\newcommand{\uc}{\mathbb{S}^1}
\newcommand{\pr}{\P_r}

\newcommand{\pf}[1]{\F{\P}{\P #1}}
\newcommand{\pft}[2]{\F{\P #1}{\P #2}}
\newcommand{\ftg}{F_t(\ga_0)}

\newcommand{\mR}{\mathbb{R}}

\def\XXint#1#2#3{{\setbox0=\hbox{$#1{#2#3}{\int}$}
    \vcenter{\hbox{$#2#3$}}\kern-.5\wd0}}

\makeatletter
\def\@citestyle{\m@th\upshape\mdseries}
\def\citeform#1{{\bfseries#1}}
\def\@cite#1#2{{%
  \@citestyle[\citeform{#1}\if@tempswa, #2\fi]}}
\@ifundefined{cite }{%
  \expandafter\let\csname cite \endcsname\cite
  \edef\cite{\@nx\protect\@xp\@nx\csname cite \endcsname}%
}{}
\makeatother
\begin{document}
\title{Curve Shortening Flows in Warped Product Manifolds}
\author{Hengyu Zhou}
\address{Department of Mathematics, Sun Yat-sen University, No. 135, Xingang Xi Road, Guangzhou, 510275, P. R. China}
\email{hyuzhou84@yahoo.com}

\begin{abstract}
   We study curve shortening flows in two types of warped product manifolds. These manifolds are $\uc \times N$ with two types of warped metrics where $\uc$ is the unit circle in $\R^2$ and $N$ is a closed Riemannian manifold. If the initial curve is a graph over $\uc$, then its curve shortening flow exists for all times and finally converges to a geodesic closed curve.
\end{abstract}
\date{\today}
\maketitle

\section{Introduction}
\subsection{Background}
  The curve shortening flow is an intriguing topic in the study of mean curvature flows. It is well-known that some properties of curve shortening flows are unique and have no correspondences for mean curvature flows of surfaces in $\R^3$. For instance, see the work of Gage-Hamilton \cite{GH86}, Grayson \cite{Gra87},\cite{Gra89b},\cite{Gra89} etc. On the other hand, curve shortening flows can be also viewed as one-dimensional mean curvature flows. An interesting example comes from the comparisons between the results of Wang \cite{WMT02} and Ma-Chen \cite{MC07}. The former work described the graphical mean curvature flows in product manifolds with dimension greater than one. The latter one studied the curve shortening flows for a ramp curve\footnote{It is a graphical curve over $\uc$ in the product manifold $\uc\times N$.} in product manifolds, which is an one-dimensional case for the former one.\\
    \indent In recent years, geometric properties of warped product manifolds are investigated by many authors. For instance, Montiel \cite{Mon99}, Brendle \cite{Bre13} and Alias-Dajczer \cite{AD07} etc. At the same time, there are growing interests in the interaction between geometric flows and warped product manifolds. For example, see Borisenko-Miquel \cite{BM12}, Tran \cite{TH15}. \\
    \indent The main motivation here is to explore the connection between warped product manifolds and curve shortening flows.
    \subsection{Main Results}
\indent Now we describe main concepts we will work with. Let $\uc$ denote the conical circle with the induced metric $d r^2$. The curve shortening flow
 $F_{t}(\ga_0):\uc \times [0, t_0)\rightarrow M$  is a smooth solution of the following quasilinear equation:
\be\label{def:pmcf}
   \F{d \ga}{dt}= \vec{H};\qquad
    \ga(.,0)=\ga_0
\ene
where $\vec{H}$ is the mean curvature vector of $F_{t}(\ga_0)$ in the {\RM} $M$ and $\ga_0$ is a closed smooth curve. \\
\indent Throughout this note, $(N,g)$ is a fixed closed {\RM} with a Riemannian metric $g$ and $(\uc, dr^2)$ is the unit circle. \\
\indent  We deal with the following two types of warped product manifolds.
\begin{Def} \label{def:wp} We call $M_{L}$ a left {\Wp} if it is $\uc \times N$ with the metric
$
 g_L= g+ \psi^2(x)dr^2
$
where $\psi:N\rightarrow \R$ is a smooth positive function.\\
 \indent We call $M_{R}$ a right {\Wp} if it is $\uc\times N$ with the metric
$
 g_R= \phi^2(r) g+ dr^2
$
where $\phi:\uc\rightarrow \R$ is a smooth positive function.
\end{Def}
\br Here we only require that $\psi$ and $\phi$ are smooth and positive. Since $S^1$ and $N$ are compact, these assumptions are sufficient for the preservation of the graphical property along curve shortening flows in our setting (see Lemma \ref{lm:estML} and Lemma \ref{lm:esMR}).
\er
The main objectives we are interested in are graphical curves and their angle functions.
\begin{Def}
\indent Suppose $\ga$ is a closed smooth curve in $M_R$ or $M_L$. Let $T$ denote a unit tangent vector of $\ga$.
\begin{enumerate}
\item In $M_L$, $\Ta_L =\la T, \pr \ra_{g_L} $ is called as the angle function of $\ga$.
\item In $M_R$, $\Ta_R= \la T, \pr \ra_{g_R} $ is called as the angle function of $\ga$.
\end{enumerate}
In addition, $\ga$ is called a graphical curve in $M_L$ or $M_R$ if it is parameterized by $(r, f(r))$ in $\uc \times N$ where $r\in \uc$ and $f:\uc \rightarrow N$ is a smooth map.
\end{Def}
 We have the following important observation. If $\ga_0$ is a graphical curve in $M_L$ or $M_R$, then its corresponding angle function $\Ta_L$ or $\Ta_R$ never vanishes respectively. Thus we can always assume $\Ta_L$ and $\Ta_R$ are positive for graphical curves.\\
\indent The main results are stated as follows.
\bt
\label{thm:MT}
Let $\ga_0$ be a graphical curve in a left warped product manifold $M_L$. Then its curve shortening flow $F_t(\ga_0)$ in $M_L$ exists for all time and converges smoothly to a totally geodesic curve.
\et
\bt\label{thm:MT2}
 Let $\ga_0$ be a graphical curve in a right warped product manifold $M_R$. Then its curve shortening flow $F_t(\ga_0)$ in $M_R$ exists for all time and converges smoothly to a totally geodesic curve.
\et
\br In fact, our proof in two results above also works for the initial conditions that $\Ta_L>0$ or $\Ta_R>0$. In this sense, we generalize Ma-Chen's result \cite{MC07} into the warped product manifolds. See Remark \ref{rm:ML}.
\er
\subsection{The plan of this note}
 Without confusion, we use $\Ta$ to denote $\Ta_L$ in $M_L$ or $\Ta_R$ in $M_R$. The process to establish Theorem \ref{thm:MT} and Theorem \ref{thm:MT2} is summarized as follows. First we establish the evolution equation of $\Ta$ along curve shortening flows in Lemma \ref{lelu}, Lemma \ref{lm:evMR}. According to the evolution equation of $\Ta$, we obtain the condition $\Ta>0$ will be preserved if the curve shortening flows exists and
    $
      \F{d}{dt}\Ta \geq \Delta \Ta +\F{1}{2}|A|^2\Ta -C
    $
    which is explained in Lemma \ref{lm:estML} and Lemma \ref{lm:esMR}.
Finally following Wang's ideas in \cite{WMT02} we exclude the finite singularity of curve shortening flows under the setting in Theorem \ref{thm:MT} and Theorem \ref{thm:MT2}. The corresponding convergence part is from Theorem D of Ma-Chen \cite{MC07}. It states that on any Riemannian manifold if the curve shortening flow for closed curves exists for all times and is not homotopic to a point, then it will converge smoothly to a totally geodesic curve. \\
  \indent This note is organized as follows. In Section 2 we record some notation and preliminary results. In Section 3 we study the evolution equation of $\Ta_L$ along curve shortening flows in $M_L$. In Section \ref{sec:thm:MT} we prove Theorem \ref{thm:MT}. In Section 5 we establish the evolution of $\Ta_R$ along curve shortening flows in $M_R$. In Section 6, we explain some minor modifications to obtain Theorem \ref{thm:MT2} based on Section 4.
\subsection{Acknowledgement} I wish to thank my advisors Prof. Zheng Huang and Prof. Yunping Jiang in Graduate Center, CUNY for their warmly discussions and encouragements during this project. I 
would also like to express my sincere gratitude to the referee for numerous corrections, many constructive suggestions and helpful comments.
 %I am grateful to Professor Jianggong You and Professor Gaofei Zhang for their interests and supports in this project.  This project is supported by
\section{One-dimensional Riemannian manifolds}
We list some notation for latter use.
Assume $M$ is a {\RM} and $\ga$ is a closed smooth curve in $M$. Suppose $F_t(\ga)$ is the curve shortening flow of $\ga$ in \eqref{def:pmcf}. Let $\{r\}$ be a local coordinate of $\ga$. Then the unit tangent vector $T$ of $F_{t}(\ga)$ can be written  as
\be \label{DefT}
   T=\F{(F_t)_*(\pf{r})}{\la(F_t)_*(\pf{r}),(F_t)_*(\pf{r})\ra^{\F{1}{2}}}
\ene
\indent
  Notice that $T$ is the only base in the tangent space of $F_{t}(\ga)$. The following facts are easily verified.
 \begin{enumerate}
  \item The mean curvature vector of $F_{t}(\ga)$ is given by $\vec{H}=\bnb_{T}T$ where $\bnb$ is the {\Cd} of $M$;
  \item Assume $\{ e_{\alpha}\}$ is a basis in the normal bundle of $F_{t}(\ga)$. Then the second fundamental form is $A(T, T)= \la \bnb_{T}T, e_\alpha\ra e_\alpha $. In addition
      \be\label{relat} |A|^2=\la \bnb_{T}T, e_{\alpha}\ra^2=|\vec{H}|^2 \ene
  \item The {\Cd} of $F_t(\ga)$ with respect to the induced metric is denoted by $\nb$. The Laplace operator of $F_t(\ga)$ with respect to the induced metric is denoted by $\Delta$. Since $\nb_{T}T=0$,
 \be
 \Delta \eta = \nb_{T}\nb_{T}\eta-\nb_{\nb_{T}T}\eta=T(T(\eta));\label{defdelta}
 \ene
 for any smooth function $\eta$ on $F_t(\ga)$.
  \end{enumerate}
 \bt[Lemma 4 in Ma-Chen \cite{MC07}]\label{auTm} Suppose $\ga$ is a closed smooth curve in {\RM} $M$ and $F_t(\ga)$ is its curve shortening flow on $[0,t_0)$. On $F_{t}(\ga)$, the mean curvature vector $\vec{H}$ and its unit tangent vector $T$ satisfy
 \be
 \bnb_{\vec{H}}T-\bnb_{T}\vec{H}=|A|^2T;
 \ene
 where $\bnb$ denotes the {\Cd} of $M$.
 \et
\section{Left warped product manifold $M_L$}
In this section, we work in the left warped product manifold $M_L$. In Definition \ref{def:wp} we define $M_L$ as $ \uc\times N$ with the metric $g_L=g+\psi^2(x)dr^2$. \\
\indent First, we demonstrate a technique result in Proposition \ref{lwp}. Then we derive the evolution equation of $\Ta_L$ along curve shortening flows in $M_L$. It gives us some estimates in Lemma \ref{lelu}, which imply that graphical curves in $M_L$ stay as graphs if their curve shortening flows exist. \\
\indent We abbreviate $\la, \ra_{g_L}$ by $\la, \ra$. The covariant derivative of $M_L$ is denoted by $\bnb$.
$D(\log\psi)$ is the gradient of $\log\psi$ with respect to the Riemannian manifold $(N,g)$. $\pr$ is the abbreviation of $\F{\P}{\P r}$.
\bpo\label{lwp} Given any two tangent vector fields $X, Y$ in $M_L$, we have
\be
\la Y, \bnb_{X}\pr\ra =\la X, D(\log\psi)\ra\la Y, \pr\ra -\la X,\pr\ra\la Y, D(\log\psi)\ra
\ene
\epo
\bp Suppose $\{x_1,\cdots, x_n\}$ is a local coordinate of $(N,g)$. This generates a new coordinate $\{x_0=r ,x_1,\cdots, x_n\}$ in $M_L$. Greek indices $\alpha, \beta, \gamma, \delta$ range from 0 to n and Latin $i,j,k,l$ from 1 to n. Define
\begin{align}
  g_{ij}&=\la \pf{x_i},\pf{x_j}\ra_g,\quad
  g_{L,\alpha\beta}=\la\pf{x_{\alpha}},\pf{x_\beta}\ra;\\
  (g^{ij})&=(g_{ij})^{-1},\quad (g_{L}^{\alpha\beta})=(g_{L,\alpha\beta})^{-1};
\end{align}
  Hence $D(\log\psi)=g^{kl}\pft{\log\psi}{x_l}\pf{x_k}$.
An obvious fact is $g_{L}^{kl}=g^{kl}, g_{L, kl}=g_{kl}$ because $g_{L}=g+\psi^2 dr^2$. Direct computations yield that the Christoffell symbols take the form $
\Gamma^{0}_{00}=0$, $\Gamma^{i}_{00}=-g^{ik}\psi\pft{\psi}{x_k}$, $
\Gamma^{i}_{0k}=0$ and $\Gamma^{0}_{0k}=\pft{\log\psi}{x_k}$.
Hence we get that
\be\label{c1}
\bnb_{\pr}\pr=-\psi  D\psi,\quad \bnb_{\pr}\pf{x_i}=\bnb_{\pf{x_i}}\pr=\pft{\log\psi}{x_i}\pr;
\ene
 A vector field $X$ can be written as
\be\label{c2}
   X=\la X, \pf{x_l}\ra g^{kl}\pf{x_k}+\la X,\pr\ra \F{1}{\psi^2}\pr;
\ene
   Combining \eqref{c1} with \eqref{c2}, we get
 \begin{align*}
 \la \bnb_{X}\pr, Y\ra &=\la \la X, \pf{x_l}\ra g^{kl}\bnb_{\pf{x_k}}\pr+\la X,\pr\ra \F{1}{\psi^2}\bnb_{\pr}\pr, Y\ra;\\
     &=\la \la X, g^{lk}\pft{\log\psi}{x_k}\pf{x_l}\ra \pr-\F{1}{\psi}\la X,\pr\ra D\psi, Y\ra;\\
     &=\la X, D(\log\psi)\ra\la \pr,Y\ra-\la X,\pr\ra\la D(\log\psi),Y\ra.
 \end{align*}
 We complete the proof.
\ep

\bl\label{lelu} Let $\ga_0$ be a closed and smooth curve in $M_L$. Suppose $F_{t}(\ga_0)$ is the curve shortening flow of $\ga_0$ on time interval $[0, t_0)$. Then the angle function $\Ta_L$ along $F_{t}(\ga_0)$ satisfies
\begin{align*}
\F{d}{dt}\Ta_L =\Delta \Ta_L +|A|^2\Ta_L +2\la \vec{H}, D(\log\psi)\ra\Ta_L-2\la \nb\Ta_L, T\ra \la T, D(\log\psi)\ra
\end{align*}
Here $\Delta$, $\nb$ denote the {\Lp} and the covariant derivative of $F_{t}(\ga_0)$ respectively. $|A|$ is the norm of the second fundamental form.
\el
\br\label{rm:ML} A special case is that $M_L$ is a product manifold $S^1\times N$ and the angle function $\Ta_L$ satisfies
$ \F{d}{dt}\Ta_L =\Delta \Ta_L +|A|^2\Ta_L $.
This is indeed Proposition 13 of Ma-Chen \cite{MC07}.
\er
\bp According to Proposition \ref{lwp}, we observe that $
\la T, \bnb_{T}\pr\ra\equiv 0 $.
Recall that $\vec{H}=\bnb_{T}T$. We compute
\begin{align}
\la \nb \Ta_L, T\ra &=T\la T, \pr\ra; \notag\\
 &=\la \bnb_{T}T,\pr\ra+\la T, \bnb_{T}\pr\ra; \notag \\
 &=\la \vec{H},\pr\ra;\label{stepo}
\end{align}
Next we compute the $t$-derivative of $\Ta_L$ as follows.
\begin{align}
\F{d}{dt}\Ta_L &=\la \bnb_{\vec{H}}T,\pr\ra +\la T, \bnb_{\vec{H}}\pr\ra \notag\\
           & = \la \bnb_{T}\vec{H}, \pr\ra +|A|^2\Ta_L +\la\vec{H}, D(\log\psi)\ra \Ta_L-\la T, D(\log\psi)\ra \la \vec{H},\pr\ra \label{step1}
\end{align}
In the second identity we apply Theorem \ref{auTm} and Proposition \ref{lwp}. On the other hand,
\begin{align}
\Delta\Ta_L &=T(T\la T,\pr\ra)\quad \mbox{by \eqref{defdelta}}\notag\\
    &=T\la \vec{H}, \pr\ra + T\la T,\bnb_{T}\pr\ra \notag\\
    &=\la \bnb_{T}\vec{H},\pr\ra+\la T, D(\log\psi)\ra\la \nb\Ta_L, T\ra-\la \vec{H}, D(\log\psi)\ra \Ta_L\label{step2}
\end{align}
where in the last step we apply Proposition \ref{lwp}.
The lemma follows from combining \eqref{stepo}, \eqref{step1} with \eqref{step2}.
\ep
The evolution equation of $\Ta_L$ is explored in the following estimation.
\bl \label{lm:estML} Assume $\Theta_L>0$ on a graphical curve $\gamma_0$ in $M_L$. Suppose its curve shortening flow $F_{t}(\ga_0)$ exists on $[0, t_0)$ where $t_0<\infty$. Then \begin{enumerate}
\item For any $t\in [0, t_0)$, $F_{t}(\ga_0)$ is graphical and
  \be\label{eq:conclusion}
   \Ta_L\geq e^{-C_L t}\min_{p\in \ga_0}\Ta_L(p)>0
  \ene
  where $C_L=\max_{q\in N}|D(\log\psi)(q)|^2$ and $|D(\log\psi)|$is the norm  of $D(\log\psi)$ with respect to $(N,g)$.
\item For $t \in [0, t_0)$, $\Ta_L$ satisfies
 \be\label{eq:eTaL}
  \F{d}{dt}\Ta_L\geq \Delta\Ta_L +\F{|A|^2}{2}\Ta_L-C_L(t_0,\ga_0);
 \ene
 where $C_L(t_0, \ga_0)=4C_L(1+\max_{q\in N}|\psi^2(q)|\F{e^{C_Lt_0}}{\min_{p\in \ga_0}\Ta_L(p)})$.
\end{enumerate}
\el
\bp To prove (1), we only have to demonstrate \eqref{eq:conclusion}. Suppose \eqref{eq:conclusion} is true and $F_{t}(\ga_0)$ loses its graphical property at some time $t$, then $\Ta_L$ has to be vanishing. It is a contradiction by \eqref{eq:conclusion}.\\
\indent Since $\Ta_L>0$ on $\ga_0$, we can assume $\Ta_L>0$ on $F_{t}(\ga_0)$ for all $t\in (0, t_1]$ with $t_1\leq t_0$. We complete the following perfect square
   \begin{align}
   |A|^2\Ta_L+2\la \vec{H},D(\log\psi)\ra  \Ta_L&\geq |A|^2\Ta_L-2|A||D\log\psi|\Ta_L;\notag\\
                   & \geq (|A|-|D(\log\psi)|)^2\Ta_L-|D(\log\psi)|^2\Ta_L \notag\\
                   & \geq - C_L\Ta_L;\label{eq:CL}
   \end{align}
   where $C_L=\max_{q\in N}|D(\log\psi)(q)|^2$. Here we also apply the estimate
   $$
    2\la \vec{H},D(\log\psi)\ra \geq - 2|\vec{H}||D(\log\psi)|=-2|A||D(\log\psi)|;
   $$ Combining Lemma \ref{lelu} with \eqref{eq:CL}, we get for all $t\leq t_1$
   \be
   \F{d}{dt}\Ta_L\geq \Delta\Ta_L - C_L\Ta_L-2\la \nb\Ta_L, T\ra\la T, D\log\psi\ra
   \ene
   The maximum principle (Theorem 4.4 in \cite{CK04}) implies that
   for $t\in [0, t_1)$ $$\Ta_L\geq e^{-C_Lt}\min_{p\in \ga_0}\Ta_L(p)\geq e^{-C_Lt_0}\min_{p\in \ga_0}\Ta_L(p)> 0$$ Repeating the above process and letting $t_1=t_0$, we establish (1). \\
   \indent In \eqref{stepo} we prove $\la \nb\Ta_L, T\ra=\la \vec{H}, \pr\ra$. Using the Cauchy inequality, we compute
   \begin{align}
     2\la \nb\Ta_L, T\ra \la T, D(\log\psi)\ra &=2\la \vec{H}, \pr\ra \la T, D(\log\psi)\ra;\\
     &\geq -2|A||D(\log\psi)||\psi|\notag\\
                        &\geq -\F{|A|^2}{4}\Ta_L-4\F{|D(\log\psi)|^2|\psi|^2}{\Ta_L};\label{eqs:3}
   \end{align}
   Here we use $\la\pr,\pr\ra=\psi^2$ and $|A|=|\vec{H}|$. On the other hand,
   \begin{align}
     2\la \vec{H}, D(\log\psi)\ra\Ta_L &\geq -2|A||D(\log\psi)|\Ta_L \notag\\
       &\geq  -\F{|A|^2}{4}\Ta_L-4|D(\log\psi)|^2\Ta_L\label{eqs:4}
       \end{align}
    With \eqref{eqs:3} ,\eqref{eqs:4} and Lemma \ref{lelu}, we obtain
    \be
      \F{d}{dt}\Ta_L\geq \Delta\Ta_L +\F{|A|^2}{2}\Ta_L -4|D(\log\psi)|^2\Ta_L-4\F{|D(\log\psi)|^2|\psi|^2}{\Ta_L};
    \ene
   According to (1), for any $t\in [0, t_0)$ we have $1\geq \Ta_L \geq e^{-C_Lt_0}\min_{p\in \ga_0}\Ta_L(p)$. Hence (2) follows from that
    $$
    -4|D(\log\psi)|^2\Ta_L-4\F{|D(\log\psi)|^2|\psi|^2}{\Ta_L}\geq -4|D(\log\psi)|^2(1+\F{|\psi|^2e^{C_Lt_0}}{\min_{p\in \ga_0}\Ta_L(p)})
    $$
    We complete the proof.
\ep
\section{The proof of Theorem \ref{thm:MT}}\label{sec:thm:MT}
 There are at least two ways to derive long time existence of mean curvature flows. The first way is to estimate directly the upper bound of the second fundamental form in terms of the evolution equation of geometric quantities (angle functions). See \cite{EH91a}, \cite{BM12} and \cite{MC07} etc. Another way is to prove the Gaussian density of any point along mean curvature flows in space-time space is $1$. White \cite{Whi05} showed that such mean curvature flows exists all times. \\
 \indent  Unfortunately, it seems that the first method fails in our case. In Wang's work \cite{WMT02} to derive the regularity of mean curvature flows, a key ingredient is to apply certain evolution equation of $*\Omega$ and demonstrate the Gaussian density\footnote{Also see Section 4 in \cite{WMT02} and Proposition 5.2 in \cite{WMT01}.} of corresponding mean curvature flows is always 1. After comparisons, we find that $*\Omega$ in \cite{WMT02} plays similar roles as that of $\Ta_L$ in $M_L$. Therefore we choose to follow the method in \cite{WMT02}.
 \subsection{The Gaussian density} First we recall some facts about Gaussian density.
 %For latter applications, we only state them in terms of curve shortening flows. For more information, we refer to \cite{WMT08} and \cite{WMT01}. \\
     We isometrically embedded  $M_L$ into a Euclidean space $\R^{n_0}$ for sufficiently large $n_0$. Suppose $\ga_0$ is a graphical curve in Theorem \ref{thm:MT}, and $F_{t}(\ga_0)$ is the curve shortening flow of $\ga_0$ existing smoothly on $[0, t_0)$. Without confusion, let $F_{t}(\ga_0)$ also denote its coordinate function in $\R^{n_0}$. The curve shortening flow equation in terms of $F_{t}(\ga_0)$ becomes
     \be
     \F{d}{dt} F_{t}(\ga_0)=\vec{H}=\tilde{H}+E;
     \ene
     where $\vec{H}$ is the mean curvature vector of $F_{t}(\ga_0)$ in $M_L$ and $\tilde{H}$ is the mean curvature vector of $F_{t}(\ga_0)$ in $\R^{n_0}$. As for $E$,  we have 
     $$
      E=-(\bnb_T T)^{\bot}=\vec{H}-\tilde{H}
     $$ 
     where $T$ is the unit tangent vector of $F_t(\gamma_0)$ and $\bot$ denote the projection into the normal bundle of $M$ (not $F_t(\gamma_0)$). Provided $M$ is compact and $T$ is a unit vector, $E$ should always be uniformly bounded independent of the position of $F_t(\gamma_0)$.\\
   \indent   A main tool to detect a possible singularity at $(y_0, t_0)$ is the backward heat kernel $\rho_{y_0,t_0}$ at $(y_0, t_0)$, which is proposed by Huisken \cite{Hui90} as follows.
    \be
    \rho_{y_0,t_0}(y,t)=\F{1}{\sqrt{4\pi(t_0-t)}}\exp(-\F{|y-y_0|}{4(t_0-t)});
    \ene
    Let $\rho_{y_0,t_0}$ be the abbreviation of $\rho_{y_0,t_0}(F_{t}(\ga_0), t)$. With direct computations, along $F_{t}(\ga_0)$ $\rho_{y_0, t_0}$ satisfies that
    \begin{align}
    \F{d\rho_{y_0,t_0}}{dt}&=-\Delta\rho_{y_0,t_0}-\rho_{y_0,t_0}(\F{|(F_t(\ga_0)-y_0)^{\bot}|^2}{4(t-t_0)^2}\notag\\
    &+\F{\la F_{t}(\ga_0)-y_0,\tilde{H}\ra}{(t_0-t)}+\F{\la (F_{t}(\ga_0)-y_0),E\ra}{2(t-t_0)});\label{eqv:rho}
    \end{align}
    where $(F_t(\ga_0)-y_0)^{\bot}$ is the normal component of $F_t(\ga_0)-y_0$ in the normal bundle of $F_t(\ga_0)$,
    and $\Delta$ is the Laplace operator of $F_{t}(\ga_0)$. \\
    \indent By the results of Huisken \cite{Hui90} and White \cite{Whi05}, for the curve shortening flow $F_{t}(\ga_0)$
    $$
    \lim_{t\rightarrow t_0, t<t_0}\int_{\ftg}\rho_{y_0,t_0} d\mu_t
    $$
    exists and is finite. The above limit is called as the \emph{Gaussian density} at $(y_0, t_0)$. Here $d\mu_t$ is the length element of $F_{t}(\ga_0)$. White \cite{Whi05} obtained that if this limit is 1, the curve shortening flow exists smoothly in a neighborhood of $(y_0, t_0)$ in the space-time space $\R^{n_0}\times [0,\infty)$.
    \subsection{The proof of Theorem \ref{thm:MT}}
    Our proof is divided into three steps.
    \bl\label{lm:estintegral} Assume $\ga_0$ is a graphical curve in $M_L$ and $F_{t}(\ga_0)$ exists smoothly on $[0,t_0)$ for $t_0<\infty$. Then
    $$
    \int_{0}^{t_0}\int_{\ftg}|A|^2\rho_{y_0, t_0}d\mu_t dt < \infty
    $$
    where $A$ is the second fundamental form of $F_{t}(\ga_0)$ in $M_L$.
    \el
    \bp Recall that
      $$
      \F{d}{dt}d\mu_t =-|\vec{H}|^2d\mu_t=-\la \tilde{H},\tilde{H}+E\ra d\mu_t;
      $$
      and $ 1\geq \Ta_L \geq e^{-C_Lt_0}\min_{p\in\ga_0}\Ta_L(p)=C_0>0$ by Lemma \ref{lm:estML}. With \eqref{eqv:rho} and (2) in Lemma \ref{lm:estML}, we compute
      \begin{align}
      &\F{d}{dt}((1-\Ta_L)\rho_{y_0,t_0}d\mu_t)\label{eq:integ}\\
      &\leq (\Delta(1-\Ta_L)\rho_{y_0,t_0}-(1-\Ta_L)\Delta\rho_{y_0,t_0})d\mu_t-C_0\F{|A|^2}{2}\rho_{y_0,t_0}d\mu_t\notag \\
      &+C_{L}(t_0,\ga_0)\rho_{y_0,t_0}d\mu_t\notag\\
      &-(1-\Ta_L)\rho_{y_0,t_0}(\F{|F^{\bot}_t(\ga_0)-y_0|^2}{4(t-t_0)^2}
   +\F{\la F_{t_0}(\ga_0)-y_0,\tilde{H}\ra}{(t_0-t)}\label{eq:nstep1}\\
   &+\F{\la (F_{t}(\ga_0)-y_0),E\ra}{2(t-t_0)}+|\tilde{H}|^2+\la \tilde{H},E\ra)\label{eq:nstep2}
      \end{align}
      We complete the perfect square in the last term of the equation above and obtain
      $$
      \eqref{eq:nstep1}+\eqref{eq:nstep2}=-(1-\Ta_L)\rho_{y_0,t_0}
      |\F{F^{\bot}_t(\ga_0)-y_0}{2(t-t_0)}+\tilde{H}+\F{E}{2}|^2+(1-\Ta_L)\rho_{y_0,t_0}\F{|E|^2}{4};
      $$
      Since $M_L$ is compact in $\R^{n_0}$, $E$ is uniformly bounded. Notice that $0\leq(1-\Ta_L)\leq 1$ and
      $\int_{\ftg}\rho_{y_0,t_0}d\mu_t$ is finite for $t\leq t_0$. Therefore, the above equation implies that
      $$
      \int_{\ftg}   [\eqref{eq:nstep1}+\eqref{eq:nstep2}] d\mu_t\leq \int_{\ftg}(1-\Ta_L)\rho_{y_0,t_0}\F{|E|^2}{4}d\mu_t=C<\infty
      $$
      Therefore the integral of \eqref{eq:integ} becomes
      \begin{align*}
      &\F{d}{dt}\int_{\ftg}((1-\Ta_L)\rho_{y_0,t_0}d\mu_t)\leq
      \int_{\ftg}(\Delta(1-\Ta_L)\rho_{y_0,t_0}-(1-\Ta_L)\Delta\rho_{y_0,t_0})d\mu_t\\
      &-C_0\int_{\ftg}\F{|A|^2}{2}\rho_{y_0,t_0}d\mu_t+C_{L}(t_0,\ga_0)\int_{\ftg}\rho_{y_0,t_0}d\mu_t
      +C
      \end{align*}
     By integrating by part the first term above vanishes. It is not hard to verify following facts for all $t\leq t_0$: (1) $\int_{\ftg}\rho_{y_0,t_0}d\mu_t$ and $\int_{\ftg}d\mu_t$ are finite; (2) $\int_{\ftg}(1-\Ta_L)\rho_{y_0,t_0}d\mu_t$ is finite. We take the integral with respect to $t$ and complete the proof.
    \ep
       We denote $F^{\lambda}_{s}(\ga_0)=\lambda(F_{t}(\ga_0)-y_0)$ for $-s=\lambda^2(t-t_0)$ where $-s\in (-\lambda^2t_0,0]$.
       \bl\label{lm:lambdajsj}
       Assume $\ga_0$ is a graphical curve in $M_L$ and $F_{t}(\ga_0)$ exists smoothly on $[0,t_0)$. Then there exists a sequence $\{s_j,\lambda_j\}_{j=1}^{\infty}$ such that $s_j\in [1,2]$, $\lim_{j\rightarrow \infty}\lambda_j=\infty$ and for any compact set $K\in \R^{n_0}$,
       \be
          \lim_{j\rightarrow \infty}\int_{F^{\lambda_j}_{s_j}(\ga_0)\cap K}|\tilde{A}_j|^2d\mu_{s_j}^{\lambda_j}=0
       \ene
       where $\tilde{A}_j$ is the second fundamental form of $F^{\lambda_j}_{s_j}(\ga_0)$ in $\R^{n_0}$, and $d\mu_{s_j}^{\lambda_j}$ is the length element of $F^{\lambda_j}_{s_j}(\ga_0)$.
       \el
       \bp According to Lemma \ref{lm:estintegral}, there exists a sequence $\{\lambda_j\}$ going to infinity such that
         \be\label{eq:limit}
         \lim_{j\rightarrow \infty}\int^{t_0-\F{1}{\lambda_j^2}}_{t_0-\F{2}{\lambda^2_j}}|A|^2\rho_{y_0,t_0}d\mu_tdt=0
         \ene
It is not hard to see that for any $p\in \ga_0$, $\lambda, s$
\begin{align}\label{sff}
\lambda A(F_{s}^{\lambda}(p))&=A(F_{t}(p)), \quad \lambda\tilde{A}(F_{s}^{\lambda}(p))=\tilde{A}(F_{t}(p))\\
& \rho_{y_0, t_0}d\mu_{t}=\rho_{0,0}d\mu^{\lambda}_{s}
\end{align}
where $\rho_{0,0}=\rho_{0,0}(F^{\lambda}_{s}(\ga_0),s)$ and $t = t_0-\F{s}{\lambda^2}$.

 According to \eqref{eq:limit}, for each $j$ there exists $s_j$ such that
      $$\F{1}{\lambda_{j}^{2}}\int_{F_{t_{j}}(\ga_0)}|A|^{2}\rho_{y_0,q_0}d\mu_{t_j}
      = \int_{F^{\lambda_j}_{s_j}(\ga_0)}|A_{j}|^{2}\rho_{0,0}d\mu^{\lambda_{j}}_{s_{j}};$$
and
    \be\label{eq:lim}
    \lim_{j\rightarrow\infty}\int_{F^{\lambda_j}_{s_j}(\ga_0)}|A_{j}|^{2}\rho_{0,0}d\mu^{\lambda_{j}}_{s_{j}}=0.
    \ene
where $t_j=t_0 -\F{s_j}{\lambda_j^2}$.
    We analyze it more carefully.
      $$ \rho_{0,0}(F^{\lambda_j}_{s_j}(\ga_0),s_j)=\F{1}{\sqrt{4\pi s_j}}\exp\big(-\F{|F^{\lambda_{j}}_{s_j}(\ga_0)|^{2}}{4\pi s_{j}}\big)
      $$
      Fix $R>0$, let $B^{R}(0)$ be the ball centered at $0$ with radius $R$ in ${\mR}^{n_0}$. Since $s_{j}\in [1,2]$, we get
      $$
     \int_{F^{\lambda_j}_{s_j}(\ga_0)}|A_{j}|^{2}\rho_{0,0}d\mu^{\lambda_{j}}_{s_{j}}\geq
      \F{1}{\sqrt{8\pi }}\exp(-\F{R^{2}}{4\pi})
      \int_{F^{\lambda_j}_{s_j}(\ga_0)}|A_{j}|^{2}d\mu^{\lambda_{j}}_{s_{j}}
      $$
     For any fixed compact set $K\subset R^{N}$ \eqref{eq:lim} leads to
         $$
          \int_{F^{\lambda_j}_{s_j}(\ga_0)\cap K}|A_{j}|^{2}d\mu^{\lambda_{j}}_{s_{j}}\rightarrow 0
         $$
    Notice that $\tilde{A}-A$ is the component of $\tilde{A}$ in the normal bundle of the left warped product manifold $M_L$ in $\R^{n_0}$. It is uniformly bounded since $M_L$ is compact. On the other hand, by \eqref{sff} we have
    \be\label{sff2}
       |\tilde{A_{j}}-A_{j}|=\F{|\tilde{A}-A|}{\lambda_{j}}\rightarrow 0 \text{ as\,} j \rightarrow \infty;
    \ene
    With \eqref{sff2} we obtain
     $$
          \int_{F^{\lambda_j}_{s_j}(\ga_0)\cap K}|\tilde{A}_{j}|^{2}d\mu^{\lambda_{j}}_{s_{j}}\rightarrow 0
      $$
         We complete the proof.
         \ep

      \bl\label{lm:rge} Assume $\ga_0$ is a graphical curve in $M_L$ and $F_{t}(\ga_0)$ is its curve shortening flow.  Let $t_0>0$ be any finite time. Then for any $(y_0, t_0)$,
        \be
         \lim_{t\rightarrow t_0, t< t_0}\int_{F_t(\ga_0)}\rho_{y_0, t_0}d\mu_t =1
        \ene
        That is, the curve shortening flow $F_t(\ga_0)$ exists for all time.
      \el
         \bp We take  $\lambda_j$ and $s_j$ in Lemma \ref{lm:lambdajsj}.
          \begin{align}
           \lim_{t\rightarrow t_0}\int_{F_t(\ga_0)}\rho_{y_0, t_0}d\mu_t &=\lim_{j\rightarrow \infty}\int_{F_{t_0-\F{s_j}{\lambda^2_{j}}}(\ga_0)}\rho_{y_0, t_0}d\mu_{t_0-\F{s_j}{\lambda^2_{j}}}\notag\\
           &=\lim_{j\rightarrow \infty}\int_{F_{s_j}^{\lambda_j}(\ga_0)}\rho_{0,0}d\mu_{s_j}^{\lambda_j}\label{eq:laststep}
            \end{align}
          We may assume that the origin is a limit point of $F_{s_j}^{\lambda_j}(\ga_0)$ otherwise the limit above is 1 and nothing needs to prove. \\
       \indent Recall that in Lemma \ref{lm:estML} we conclude that $\Ta_L\geq C >0$ on $[0, t_0)$. On the other hand, $F_t(\ga_0)$ is the graph of $f_t$ in $M_L$. And $\Ta_L$ is written as
        \be\label{eq:keystep:1}
        \Ta_L = \F{1}{\sqrt{1+\psi^2(x)|df_t|^2}}
        \ene
  Therefore we conclude that \emph{$df_t$ is uniformly bounded on $[0, t_0)$}.  \\
 \indent Denote $f_{t-\F{s_j}{\lambda_j}}$ by $f_j$. Therefore $F_{s_j}^{\lambda_j}(\ga)$ is the graph of $\tilde{f}_j=\lambda_j f_j$ which is defined on $\lambda_j \uc\subset \R^{n_0}$.  It is not hard to see that $d\tilde{f}_j$ is also uniformly bounded. Now our assumptions on $F_{s_j}^{\lambda_j}$ imply $\lim_{j\rightarrow 0}\tilde{f}_j(0)=0$. Therefore we may assume $\tilde{f}_j\rightarrow \tilde{f}_\infty$ in $C^\alpha$ on compact sets, where $\tilde{f}_\infty$ is a map on $\R^1$.
 On the other hand, with similar computations as in Equation (29) on Page 30 of \cite{Ilm95}, we have
         \be \label{sob 2}
                  |\tilde{A}_{j}|\leq |\nabla d \tilde{f}_{j}|\leq (1+|d\tilde{f}_{j}|^{2})^{\F{3}{2}}|\tilde{A}_{j}|;
         \ene
 Here $\nabla$ is the covariant derivative of $F_{s_j}^{\lambda_j}$. From Lemma \ref{lm:lambdajsj}, $\tilde{f}_j\rightarrow \tilde{f}_\infty$ in $C^{\alpha}\cap W^{1,2}_{loc}$ and the second derivative of  $\tilde{f}_\infty$ is vanishing. Then $f_{\infty}$ is an one-dimensional linear map.
 This implies that $F^{\lambda_j}_{s_j}(\ga_0)\rightarrow F_{-1}^{\infty}$ as Radon measures and $F_{-1}^{\infty}$ is the graph of a linear map. Finally
  $$
  \lim_{j\rightarrow \infty}\int_{F_{s_j}^{\lambda_j}(\ga_0)}\rho_{0,0}d\mu_{s_j}^{\lambda_j}=\int_{F_{-1}^{\infty}}\rho_{0,0}d\mu_{-1}^\infty=1
  $$
  Our lemma follows from \eqref{eq:laststep}.
\ep
 Now we are ready to conclude Theorem \ref{thm:MT}.
\bp According to Lemma \ref{lm:rge} and White's regularity results \cite{Whi05} $(y_0, t_0)$ is a regular point of the curve shortening flows $F_{t}(\ga_0)$ in $M_L$. Therefore, for a graphical $\ga_0$ the curve shortening flow $F_{t}(\ga_0)$ in $M_L$ exists for all time. Theorem \ref{thm:MT} follows from Theorem D in (Ma-Chen,\cite{MC07}). It says that as long as a curve shortening flow exists for all time, it will converge smoothly to a geodesic curve.
\ep
%\br\label{rm:logre}  In the proof of Theorem \ref{thm:MT}, the only fact we need is Lemma \ref{lelu}. We summarize their logic structure as follows:
 % \begin{enumerate}
  %\item [(1)]Lemma \ref{lm:estML} implies Lemma \ref{lm:estintegral}.
  %\item [(2)]Lemma \ref{lm:lambdajsj} follows from Lemma \ref{lm:estML}.
  % \item [(3)] In the final step, a key ingredient is equation \eqref{eq:keystep:1} which is from (1) in Lemma \ref{lm:estML}.
       %\end{enumerate}
%\er
\section{Right warped product manifold $M_R$}\label{sec:thm:MT2}
In this section, we study curve shortening flows in the right warped product manifold $M_R$. In Definition \ref{def:wp} we define $M_R$ as $\uc\times N$ with the metric $g_R=\phi^2(r)g + dr^2$.\\
       % \indent The sketch of this section is similar as that of the previous section. First, we establish a technique % lemma. Then we derive the evolution equation of $\Ta_R$ along curve shortening flows. Then we will see some % % estimations which implies that $\Ta_R>0$ will be preserved along curve shortening flows. \\
\indent In this section, we denote $\la ,\ra_{g_R}$ by $\la, \ra$. Let $\bnb$ be the {\Cd} of $M_R$. $(\log\phi)', (\log\phi)''$ are the first derivative and second derivative of $\log\phi$ respectively.
\bpo\label{pro:teMR}
 Given two tangent vector fields $X, Y$ in $M_R$, we have
\be
\la Y, \bnb_{X}\pr\ra =(\log\phi)'(\la X,Y\ra -\la X,\pr\ra \la Y, \pr\ra)
\ene
\epo
\bp  A well-known fact is that $\phi(r)\P_r$ is a conformal vector field (see Montiel \cite{Mon99}), i.e $\bnb_X (\phi(r)\P_r)=\phi'(r)X$ for any smooth vector field $X$. Then
 \begin{align*}
   \la Y, \bnb_X\P_r\ra&=\la Y,\bnb_X(\F{\phi(r)\P_r}{\phi(r)})\ra\\
                      &=\F{\phi'(r)}{\phi(r)}(\la X, Y\ra-\la Y,\P_r\ra\la X,\P_r\ra)
   \end{align*}
\ep
With this result we obtain the evolution equation of $\Ta_R$.
\bl\label{lm:evMR} Let $\ga_0$ be a closed smooth curve in $M_R$. Suppose $F_{t}(\ga_0)$ is the curve shortening flow of $\ga_0$ on $[0, t_0)$. Then the angle function $\Ta_R$ satisfies
\be\label{eq:evMR}
\F{d}{dt}\Ta_R =\Delta \Ta_R +|A|^2\Ta_R +2(\log\phi)'\Ta_R\la \nb\Ta_R, T\ra -(\log\phi)''\Ta_R(1-\Ta_R^2);
\ene
Here $\Delta$, $\nb$ denote the {\Lp} and the covariant derivative of $F_{t}(\ga_0)$ respectively. $|A|$ is the norm of the second fundamental form.
\el
\bp Recall that $\Ta_R =\la T,\pr\ra$. From Proposition \ref{pro:teMR}, we observe an symmetry property as follows. For any two tangent vector fields $X, Y$,
\be \label{eq:sym}
\la X, \bnb_{Y}\pr\ra = \la Y, \bnb_{X}\pr\ra
\ene
Applying Theorem \ref{auTm} and \eqref{eq:sym}, we compute $\F{d}{dt}\Ta_R$ as
\begin{align}
\F{d}{dt}\Ta_R &=\la \bnb_{\vec{H}}T,\pr\ra +\la T, \bnb_{\vec{H}}\pr\ra \notag\\
               &=\la \bnb_{T}\vec{H}, \pr\ra +|A|^2\Ta_R +\la \vec{H}, \bnb_{T}\pr\ra\label{eq:n1}
\end{align}
On the other hand, by \eqref{defdelta} we have
\begin{align}
\Delta\Ta_R &= T(T\la T, \pr\ra)\notag\\
            &=T\la \vec{H}, \pr\ra + T\la T,\bnb_{T}\pr\ra\notag\\
            &=\la \bnb_{T}\vec{H}, \pr\ra +\la \vec{H}, \bnb_{T}\pr\ra + T((\log\phi)'(1-\Ta_R^2))\notag\\
            &=\la \bnb_{T}\vec{H}, \pr\ra +\la \vec{H}, \bnb_{T}\pr\ra+(\log\phi)''\Ta_R(1-\Ta_R^2)\notag\\
            &-2(\log\phi)'\Ta_R\la \nb\Ta_R, T\ra;\label{eq:n2}
\end{align}
because $\la T, \bnb_{T}\pr\ra=(\log\phi)'(1-\Ta_R^2)$  by Proposition \ref{pro:teMR}. The Lemma follows from combining \eqref{eq:n1} with \eqref{eq:n2}.
\ep
We have the following estimations for $\Ta_R$.
\bl\label{lm:esMR} Assume $\Theta_R>0$ on a graphical curve $\gamma_0$ in $M_R$. If its curve shortening flow $F_{t}(\ga_0)$ exists on $[0, t_0)$ where $t_0<\infty$,
 \begin{enumerate}
\item For any $t\in [0, t_0)$, $F_t(\ga_0)$ are graphical curves and
  \be\label{eq:cv}
  \Ta_R\geq e^{-C_R t}\min_{p\in \ga_0}\Ta_R(p)>0
  \ene
  where $C_R=\max_{r\in \uc}|(\log\phi)''(r)|$.
\item For $t \in [0, t_0)$, $\Ta_R$ satisfies the following inequality:
 \be\label{eq:estTaL}
  \F{d}{dt}\Ta_R\geq \Delta\Ta_R +\F{|A|^2}{2}\Ta_R-C_R(\phi);
 \ene
 where $C_R(\phi)=\max_{r\in \uc}\{4((\log\phi)')^2+|(\log\phi)''|\}(r)$.
\end{enumerate}
\el
\bp With the same reason as in Lemma \ref{lm:estML}, it is only sufficient to prove \eqref{eq:cv} to conclude (1). \\
 \indent Since $\Ta_R>0$ on $\ga_0$, we can assume $\Ta_R>0$ on $F_{t}(\ga_0)$ for some $t_1\in (0, t_0]$. It is easy to see that for $t\in (0, t_1)$,
\be\label{eq:dfa}
-(\log\phi)''\Ta_R(1-\Ta_R^2)\geq -|\log\phi''|\Ta_R\geq-C_R\Ta_R
\ene
since $0<\Ta_R\leq 1$ for $t\in [0, t_1)$. Here $C_R=\max_{r\in \uc}|(\log\phi)''(r)|$.
 According to  Lemma \ref{lm:evMR}, for $t\in [0, t_1)$ we have
   \be
   \F{d}{dt}\Ta_R\geq \Delta\Ta_R - C_R\Ta_R+2(\log\phi)'\Ta_R^2\la \nb\Ta_R, T\ra ;
   \ene
    By the maximum principle (Theorem 4.4 in \cite{CK04}) we get
    $$\Ta_R\geq e^{-C_Rt}\min_{p\in \ga_0}\Ta_R(p)\geq e^{-C_Rt_0}\min_{p\in \ga_0}\Ta_R(p)>0$$ for $t\in [0, t_1)$.
     Repeating the above process, we can choose $t_1=t_0$ and therefore obtain (1). \\
   \indent By Proposition \ref{pro:teMR}, we have
   $$
    \la \nb\Ta_R, T\ra=T\la T, \pr\ra = \la \vec{H}, \pr\ra +(\log\phi)'(1-\Ta_R^2);
   $$
       Using the Cauchy inequality on $-2(\log\phi)'\Ta_R^2\la \nb \Ta, T\ra$ and $1\geq\Ta_R>0$, we get
   \begin{align}
     -2(\log\phi)'\Ta_R^2\la \nb \Ta, T\ra &= -2(\log\phi)'\Ta_R^2\la \vec{H}, \pr\ra -2((\log\phi)')^2\Ta_R^2(1-\Ta_R^2)\notag\\
                        &\geq -2(\log\phi)'|A|\Ta_R-2((\log\phi)')^2\Ta_R\notag \\
                        &\geq -\F{|A|^2\Ta_R}{2}-4((\log\phi)')^2\label{eq:midR}
   \end{align}
   Here we use $|\vec{H}|=|A|$ since $F_{t}(\ga_0)$ are curves.
    Plugging \eqref{eq:dfa} and \eqref{eq:midR} into \eqref{eq:evMR}, we obtain
    \be
      \F{d}{dt}\Ta_R\geq \Delta\Ta_R +|A|^2\Ta_R -\F{|A|^2\Ta_R}{2}-4((\log\phi)')^2-|(\log\phi)''|;
    \ene
 Let $C_R(\phi)$ denote $\max_{r\in \uc}\{4((\log\phi)')^2+|(\log\phi)''|\}(r)$, we get (2). We complete the proof.
\ep
\section{The proof of Theorem \ref{thm:MT2}}
The proof of Theorem \ref{thm:MT2} is very similar as that of Theorem \ref{thm:MT} with minor modifications. We list the key ingredients in the proof. Their derivations are skipped and can be easily written according to Section 4.
  First, Lemma \ref{lm:esMR} will imply
  \bl \label{lm:estMR} Assume $\ga_0$ is a graphical curve in $M_R$ and $F_{t}(\ga_0)$ exists smoothly on $[0,t_0)$ for $t_0<\infty$. Then
    $$
    \int_{0}^{t_0}\int_{\ftg}|A|^2\rho_{y_0, t_0}d\mu_t dt < \infty
    $$
    where $A$ is the second fundamental form of $F_{t}(\ga_0)$ in $M_R$.
\el
 Then Lemma 6.1 indicates the following result.
    \bl
       Assume $\ga_0$ is a graphical curve in $M_R$ and $F_{t}(\ga_0)$ exists smoothly on $[0,t_0)$. Then there exists $\{s_j,\lambda_j\}_{j=1}^{\infty}$ such that $s_j\in [1,2]$, $\lim_{j\rightarrow \infty}\lambda_j=\infty$ and for any compact set $K\in \R^{n_0}$,
       \be
          \lim_{j\rightarrow \infty}\int_{F^{\lambda_j}_{s_j}(\ga_0)\cap K}|\tilde{A}_j|^2d\mu_{s_j}^{\lambda_j}=0
       \ene
       where $\tilde{A}_j$ is the second fundamental form of $F^{\lambda_j}_{s_j}(\ga_0)$ in $\R^{n_0}$, and $d\mu_{s_j}^{\lambda_j}$ is the length element of $F^{\lambda_j}_{s_j}(\ga_0)$.
       \el
       The way to choose $s_j, \lambda_j$ is the same as in Section 4.
   One difference is that equation \eqref{eq:keystep:1} is replaced by
       $
       \Ta_R = \F{1}{\sqrt{1+|df_t|^2(\phi(r))^{-2}}}
       $.
       Here $|df|$ is the norm of $df$ with respect to the Riemannian manifold $(N,g)$. The conclusion that $df$ is uniformly bounded still works under the setting of Theorem \ref{thm:MT2}. Then we establish Theorem \ref{thm:MT2} with the similar derivation as in Section 4.

\bibliographystyle{abbrv}	
\bibliography{csf}
\end{document}